\documentclass[11pt]{amsart}
\usepackage{amssymb}
\usepackage{amscd}
\usepackage{curves}
\usepackage{epsfig}
\usepackage[pass]{geometry}
\usepackage{graphicx}
\usepackage{pstricks}
\usepackage{pstricks-add}
\usepackage{pst-grad}
\usepackage{pst-plot}
\usepackage{verbatim}
\usepackage{latexsym,eucal,amsfonts,amssymb,amsmath,graphicx}
\usepackage{auto-pst-pdf}
\usepackage{xcolor}

\numberwithin{equation}{section}
\newtheorem{theorem}{Theorem}[section]

\newtheorem{lemma}[theorem]{Lemma}
\newtheorem{prop}[theorem]{Proposition}

\def \bpf {\begin{proof}}
\def \epf {\end{proof}}
\def \beq {\begin{equation*}}
\def \eeq {\end{equation*}}
\def \bsp{\begin{split}}
\def \esp{\end{split}}
\def \beqq {\begin{equation}}
\def \eeqq {\end{equation}}
\def \mca {{\mathcal A}}
\def \mcb {{\mathcal B}}

\def \mce {{\mathcal E}}
\def \mcf {{\mathcal F}}

\def \mck {{\mathcal K}}

\def \mcm {{\mathcal M}}
\def \mcn {{\mathcal N}}
\def \mco {{\mathcal O}}
\def \mcp {{\mathcal P}}

\def \mcs {{\mathcal S}}

\def \mcu {{\mathcal U}}

\def \mcw {{\mathcal W}}
\def \mcx {{\mathcal X}}
\def \mcy {{\mathcal Y}}
\def \mcz {{\mathcal Z}}

\def \mbr {{\mathbb R}}

\def \comp {\operatorname{comp}}

\def \loc {\operatorname{loc}}

\def \det {\operatorname{det}}

\def \diag{\textrm{Diag}}

\def \supp {\text{supp }}

\def \eps {\epsilon}   
   
\def \La {\Lambda}

\def \p {\partial}

\def \eps {\epsilon}
\def \det {\text{det}}

\def \ha {\frac{1}{2}}

\def \fnf {\frac{n}{4}}
\def \WF {\operatorname{WF}}
\def \loc {\operatorname{loc}}
\def \comp {\operatorname{comp}}

\def \ba {\begin {eqnarray*} }
\def \ea {\end {eqnarray*} }

\begin{document}
\title{Inverse problems for  quadratic derivative nonlinear wave equations}

\author{Yiran Wang}
\address{Yiran Wang
\newline
\indent Department of Mathematics, University of Washington 
\newline
\indent Box 354350, Seattle, WA  98195-4350 
\newline
\indent and 
\newline
\indent Institute for Advanced Study, the Hong Kong University of Science and Technology
\newline
\indent Lo Ka Chung Building, Lee Shau Kee Campus, Clear Water Bay, Kowloon, Hong Kong}
\email{wangy257@math.washington.edu}
\author{Ting Zhou}
\address{Ting Zhou
\newline
\indent Department of Mathematics, Northeastern University 
\newline
\indent 360 Huntington Ave., Boston, MA 02115
}
\email{t.zhou@neu.edu}

\begin{abstract} 
For semilinear wave equations on Lorentzian manifolds with quadratic derivative nonlinear terms, we study the inverse problem of determining the background Lorentzian metric. Under some conditions on the nonlinear term, we show  that from the source-to-solution map, one can determine the Lorentzian metric up to diffeomorphisms.  
\end{abstract}

\maketitle

\section{Introduction}
Let $g$ be a time oriented Lorentzian metric  on $M = \mbr^{1 + 3}$ with signature $(-, +, +, +)$. Let $x = (x^0, x^1, x^2, x^3), x^0 = t$ be the coordinate for $M$. We can write $g = \sum_{i, j = 0}^{3} g_{ij}(x) dx^idx^j$. The Laplace-Beltrami operator is given by
\beq
\square_g = (-\det g(x))^{-\ha} \sum_{i, j = 0}^{3} \p_i ((-\det g(x))^\ha g^{ij}(x) \p_j), \ \ \p_i = \frac{\p}{\p x^i}.
\eeq
We consider the following semi-linear wave equation
\beqq\label{eqnlw}
\square_g u + w(x, u, \nabla_g u) = 0,
\eeqq
where $\nabla_g u$ denotes the gradient of $u$: $(\nabla_g u)^i = \sum_{j = 0}^3 g^{ij}\p_j u$,  $w(x, u, \xi)$   is smooth in $x, u$ and quadratic in $\xi\in \mbr^4$. This type of equation is usually called quadratic derivative nonlinear wave equations, which serves as a prototype for many equations in mathematical physics e.g.\ the wave map equations, see \cite[Chap.\ 6]{Tao}. The forward problem (local and global well-posedness) of \eqref{eqnlw} has been studied extensively in the literature, see e.g.\ \cite{Tao}. In this work, we study the inverse problem of determining the Lorentzian metric $g$ using equation \eqref{eqnlw}. Such problem has been considered for semilinear wave equations with no derivatives, see \cite{KLU, LUW}. Our first motivation is to take into account the derivative terms. Moreover, we'd like to get some insights into inverse problems for more complicated systems with similar nonlinear terms, for instance the Einstein equations in wave gauge. In fact, there are some interesting effects related to the null forms.

We introduce some notions to state the problem. For $p, q\in M$, we denote by $p\ll q$ ($p< q$) if $q$ is on a future pointing-time like (causal) curve from $p$ and $q\neq p$. We write $p\leq q$ if $p = q$ or $p<q$. We denote the chronological (causal) future of $p\in M$ by $I^+(p) = \{q\in M: p\ll q\}$ ($J^+(p) = \{q\in M: p\leq q\}$). The chronological and causal past are denoted by $I^-(p)$ and $J^-(p)$ respectively. We denote $J(p, q) = J^+(p)\cap J^-(q)$ and $I(p, q) = I^+(p)\cap I^-(q)$. For any set $A\subset M$, we let $J^\pm(A) = \cup_{p\in A}J^\pm(p)$. In this work, we assume that $g$ is globally hyperbolic. According to Bernal and S\'anchez \cite{BS}, this means that there is no closed causal paths in $M$ and for any $p, q\in M$ and $p<q$, the set $J(p, q)$ is compact. 
 
\begin{figure}[htbp]
\centering
\scalebox{0.7} 
{
\begin{pspicture}(0,-4.68)(6.8418946,4.68)
\psellipse[linewidth=0.04,dimen=outer](3.11,-0.13)(0.55,3.69)
\psbezier[linewidth=0.04](3.0557144,-4.66)(3.26,-3.04)(3.2428572,-2.26)(3.1114285,-0.18)(2.98,1.9)(3.0557144,3.58)(3.2414286,4.66)
\psellipse[linewidth=0.04,linestyle=dashed,dash=0.16cm 0.16cm,dimen=outer](3.1,-0.23)(3.1,0.61)
\psdots[dotsize=0.12](3.08,2.84)
\psdots[dotsize=0.12](3.18,-3.18)
\psline[linewidth=0.04cm,linestyle=dashed,dash=0.16cm 0.16cm](3.12,2.84)(6.1,-0.08)
\psline[linewidth=0.04cm,linestyle=dashed,dash=0.16cm 0.16cm](3.04,2.84)(0.08,-0.1)
\psline[linewidth=0.04cm,linestyle=dashed,dash=0.16cm 0.16cm](6.16,-0.3)(3.18,-3.2)
\psline[linewidth=0.04cm,linestyle=dashed,dash=0.16cm 0.16cm](0.04,-0.3)(3.16,-3.2)
\usefont{T1}{ptm}{m}{n}
\rput(4.5,2.965){$p_+ = \hat\mu(s_+)$}
\usefont{T1}{ptm}{m}{n}
\rput(4.5,-3.255){$p_- = \hat\mu(s_-)$}
\usefont{T1}{ptm}{m}{n}
\rput(1.5614551,-0.19){$I(p_-, p_+)$}
\usefont{T1}{ptm}{m}{n}
\rput(2.95,-2.235){$V$}
\usefont{T1}{ptm}{m}{n}
\rput(3.5071455,4.385){$\hat\mu$}
\end{pspicture} 
}
\caption{Illustration of the inverse problem. $\hat\mu([-1, 1])$ is a time-like geodesic.  $V$ is an open neighborhood of $\hat\mu([s_-, s_+]), -1< s_-<s_+<1$. The source $f$ is supported in $V$ and the measurement $u = L(f)$ is in $V$. The inverse problem is to determine the metric $g$ in $I(p_-, p_+)$ bounded by the dashed curves.}
\label{figinv}
\end{figure}
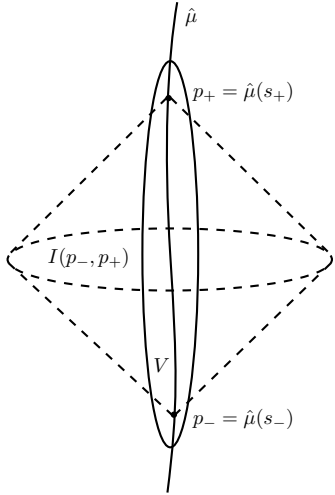

Let $\hat \mu(t) \subset M$ be a time-like geodesic where $t\in [-1, 1]$. In general relativity, $\hat \mu$ represents a freely falling observer. We consider the equation \eqref{eqnlw} near $\hat \mu$. Without loss of generality, we assume $\hat \mu(-1) \in \{0\}\times \mbr^3$. Take $T_0 >0$ such that $\hat \mu([-1, 1]) \subset M(T_0) \doteq (-\infty, T_0)\times \mbr^3$. Let $p_\pm = \hat\mu(s_\pm), -1<s_-<s_+<1$ be two points on $\hat \mu$ and $V$ be an open relatively compact neighborhood of $\hat \mu([s_-, s_+])$ with $V\subset M(T_0)$. We consider \eqref{eqnlw} with sources
\beqq\label{eqsem}
\begin{gathered}
\square_{g} u(x) + w(x, u(x), \nabla_g u(x))  = f(x), \text{ on } M(T_0),\\
u = 0 \text{ in } M(T_0)\backslash J^+_{g}(\supp(f)),
\end{gathered}
\eeqq
where $f$ is compactly supported in $V$. We discuss the well-posedness of \eqref{eqsem} in Section \ref{asym}. 
By taking a complete Riemannian metric $g^+$ on $M$, we can define (semi)norms of $C^m$ and Sobolev spaces on $M$. Then there exists a unique solution $u\in H^4(M(T_0))$ of \eqref{eqsem} for $f$ small in $C^4(M)$ and compactly supported in $V$. For such $f$, we define the source-to-solution map as 
\beq
L(f) = u|_{V}.
\eeq
We shall regard $f$ as the controllable source that generates the nonlinear wave  $u$ and $L(f)$ is the measurement. The inverse problem we study in this work is that given the source-to-solution map $L$, can one determine the metric $g$   on the set $I(p_-, p_+)$, which is the set where the waves can propagate to from $\hat \mu$ and return   to $\hat \mu$? See Fig.\ \ref{figinv}.

This type of inverse problem has been studied for semilinear wave equations of the form
\beqq\label{semieq}
\square_g u(x)+ H(x, u(x)) = f(x),
\eeqq
starting from the quadratic nonlinearity i.e.\ $H(x, u) = a(x)u(x)^2$ in Kurylev-Lassas-Uhlmann \cite{KLU}. The main result Theorem 1.5 of \cite{KLU} states that if  $a(x)$ is non-vanishing, one can determine the conformal class of the metric $g$ from the source-to-solution map. The results are extended to general nonlinear terms with no derivatives in Lassas-Uhlmann-Wang \cite{LUW}. Similar inverse problems have been studied for the Einstein equations in general relativity, see \cite{KLU1, KLU2}.

As introduced in \cite{KLU1}, a key idea to solve this type of inverse problems is to produce point source like singularities in $I(p_-, p_+)$ from the nonlinear interaction of four progressive waves. The inverse problem can be solved by observing the leading singularities from the point source. The phenomena, known as nonlinear interaction of singularities, was actively studied in the 80's and 90's, see for example Melrose-Ritter \cite{MR}, Rauch-Reed \cite{RR} and Beals \cite{Bea}. For equations with quadratic forms, thanks to the work of Klainerman-Machedon \cite{KM1, KM2}, we know that null forms have certain smoothing effects on the regularity of solutions. This causes additional difficulties as the leading singularities we expect to observe might actually vanish. 

Let\rq{}s recall that for the Minkowski space-time $(\mbr^4, h), h = -dt^2 + \sum_{j = 1}^3 (dx^j)^2$, a quadratic form $w(\xi, \eta)$ satisfies the null condition if
\beq
w(\xi, \xi) = 0  \text{ for  any $\xi \in \mbr^4$ such that $h(\xi, \xi) = 0$},
\eeq
see e.g.\ \cite[Def.\ 3.1]{So}. On Lorentzian manifold $(M, g)$, quadratic forms can be regarded as controvariant two tensor fields  i.e. $w\in C^\infty(M; T^*M\otimes T^*M)$. We call $w$ a quadratic null form if for any $x\in M$,  
\beqq\label{cond}
w(x, \xi, \xi) = 0 \text{ for  any $\xi \in T_xM$ such that $g(x, \xi, \xi) = 0$}. \tag{N}
\eeqq
Consider the nonlinear term $w$ in \eqref{eqsem}. We can write $w$ in Taylor expansion in $u$ as
\beq
w(x, u, \xi) =  \mcn_0(x, \xi) + u \mcn_1(x, \xi) +  u^2 \mcm(x,  \xi) + o(|u|^2\cdot|\xi|^2), 
\eeq
where $\mcn_0, \mcn_1, \mcm$ are quadratic forms. In this work, we make the assumption on $w$ that
\beqq\label{eqw}
\text{ $\mcn_0, \mcn_1$ are null forms and $\mcm$ is not null}. \tag{A}
\eeqq
For quadratic forms $w^{(j)}, j = 1, 2$, we denote the corresponding terms in the Taylor expansion by $\mcn_0^{(j)}, \mcn_1^{(j)}$ and $\mcm^{(j)}$. Our main result is
\begin{theorem}\label{main}
Let $g^{(j)}, j = 1, 2$ be two globally hyperbolic Lorentzian metrics on $M = \mbr^{1+3}$. Let $\hat \mu^{(j)}(t)$ be time-like geodesics on $(M, g^{(j)})$ where $t\in [-1, 1]$. Let $M(T_0) = (-\infty, T_0)\times \mbr^3, T_0 > 0$ be such that $\hat \mu^{(j)}([-1, 1])  \subset M(T_0)$. Let $V^{(j)}\subset M(T_0)$ be an open relatively compact neighborhood of $\hat \mu^{(j)}([s_-, s_+])$ where $-1<s_-<s_+<1$. Consider the nonlinear wave equations with source terms
\beq
\begin{gathered}
\square_{g^{(j)}} u(x) + w^{(j)}(x, u(x), \nabla_{g^{(j)}} u(x) )  = f(x), \text{ on } M(T_0),\\
u = 0 \text{ in } M(T_0)\backslash J^+_{g^{(j)}}(\supp(f)),
\end{gathered}
\eeq
where $\supp(f)\subset V^{(j)}$. Suppose that $w^{(j)}$ satisfy  the assumption \eqref{eqw}. If there is a diffeomorphism $\Phi: V^{(1)}\rightarrow V^{(2)}$ such that $\Phi(p^{(1)}_\pm) = p^{(2)}_\pm$ and the source-to-solution maps $L^{(j)}$ satisfy
\beq
((\Phi^{-1})^*\circ L^{(1)}\circ \Phi^*)(f) = L^{(2)}(f)
\eeq
for all $f$  in a small neighborhood of the zero function in $C_0^4(V^{(2)})$, then
\begin{enumerate}
\item there is a diffeomorphism $\Psi: I(p^{(1)}_-, p^{(1)}_+)\rightarrow I(p^{(2)}_-, p^{(2)}_+)$ such that $g^{(1)}$ is conformal to $g^{(2)}$ i.e.\ $\Psi^*g^{(2)} = e^{2\gamma} g^{(1)}$ in $I(p^{(1)}_-, p^{(1)}_+)$ with $\gamma$ smooth.  
\item In addition, if $\mcm^{(1)} = \mcm^{(2)}$ are independent of $x$ or $g^{(1)}$ and $g^{(2)}$ are Ricci flat, the diffeomorphism $\Psi$ is an isometry. 
\end{enumerate}
\end{theorem} 

We remark that when we analyze the singularities generated from the nonlinear interaction of four linear waves, it turns out that for nonlinear terms satisfying \eqref{eqw}, the leading singularities are generated from the term $u^2\mcm$ and not from the null forms, which can be regarded as a manifest of the smoothing effects. The result could potentially help us to simplify the analysis for more complicated systems  as one can factorize the null forms out of the quadratic nonlinear term. In this work, we deal with the scalar equation \eqref{eqsem} as an example. The method should be applicable to systems as well.  

We also remark that we take $M = \mbr^{1+3}$ mainly for simplicity. Actually, for any globally hyperbolic Lorentzian manifold $(M, g)$, it is proved in \cite{BS0} that $(M, g)$ can be identified with the product manifold 
\beq
\mbr\times N \text{ with metric } g = -\beta(t, y) dt^2 + \kappa(t, y),
\eeq
where $N$ is a $3$ dimensional manifold, $\beta$ is smooth and $\kappa$ is a family of Riemannian metrics on $N$ smoothly depending on $t$. The generalization of Theorem \ref{main} to $(M, g)$ is immediate.

The paper is organized as follows. In Section \ref{asym}, discuss the well-posedness of the equation \eqref{eqsem} and derive the asymptotic expansion of the solution when the source depends on a small parameter. In Section \ref{secdist}, we construct distorted plane waves which are generalizations of the traveling waves and we analyze their nonlinear interactions in Section \ref{secinter}. We finish the proof of the main theorems in Section \ref{secinv}.

\section{Local well-posedness and the asymptotic expansion}\label{asym}
First, we prove the decomposition of a quadratic null form on $(M, g)$. This is known for the Minkowski case when the quadratic form dose not depend on $x$, see \cite[Section II.3]{So}. 
\begin{lemma}\label{biform}
Let $w  \in C^\infty(M; T^*M\otimes T^*M)$ be a quadratic form satisfying condition \eqref{cond}. For any local coordinates $(x^i)_{i = 0}^3$ of $M$, we can find $C_\bullet(x)$ smooth such that $w$ can be written as
\beqq\label{eqwexpan}
w(x, \xi, \eta) = C_0(x) g(\xi, \eta) + \sum_{0\leq a< b\leq 3} C_{ab}(x) w_{ab}(\xi, \eta), \ \ \xi, \eta\in T^*M,
\eeqq
where $
w_{ab}(\xi, \eta) = \xi^a\eta^b - \xi^b\eta^a, \ \ 0\leq a < b\leq 3.$
\end{lemma}
\bpf
Let $\mcb$ be the set of quadratic forms on $\mbr^{4}$ as a vector space so $\mcb$ is a $16$-dimensional vector space. In the following, we identify quadratic forms with $4\times 4$ matrices and do not distinguish their notations. We use the following basis for $\mcb$ which is easy to verify:
\begin{enumerate}
\item $E^{\alpha\beta}, 0\leq \alpha < \beta \leq 3$ defined by $E^{\alpha\beta}_{ij} = 1$ if $i = \alpha, j = \beta$; $E^{\alpha\beta}_{ij} = -1$ if $i = \beta, j =\alpha$; otherwise $E^{\alpha\beta}_{ij} = 0$. In particular, $(E^{\alpha\beta})^T = -E^{\alpha\beta}$ is anti-symmetric where $T$ denotes the transpose of matrix. Notice that $E^{ab}$ are the same as $w_{ab}$. 

\item   $F^{\alpha\beta}, 0\leq \alpha < \beta \leq 3$ defined by $F^{\alpha\beta}_{ij} = 1$ if $i = \alpha, j = \beta$ and $i = \beta, j =\alpha$; otherwise $F^{\alpha\beta}_{ij} = 0$. In particular, $(F^{\alpha\beta})^T = F^{\alpha\beta}$ is symmetric.

\item $G^\alpha, \alpha = 0, 1, 2, 3$ defined by $G^\alpha_{ij} = 1$ if $i = j = \alpha$ and otherwise $G^{\alpha}_{ij} = 0.$ So $G^\alpha$ are diagonal.
\end{enumerate}

Let $x = (x^0, \cdots, x^3)$ be the local coordinate for $M$. For $w\in C^\infty(M; T^*M\otimes T^*M) \simeq C^\infty(M; \mcb)$, we can express $w$ using the above basis i.e.\
\beq
w(x) = \sum_{0\leq \alpha<\beta\leq 3}a_{\alpha\beta}(x)E^{\alpha\beta} + \sum_{0\leq \alpha<\beta\leq 3}b_{\alpha\beta}(x)F^{\alpha\beta} + \sum_{\alpha = 0}^3 c_\alpha(x) G^\alpha,
\eeq
where $a_\bullet, b_\bullet$ and $c_\bullet$ are all smooth functions. The metric $g$ can be identified as a $4\times 4$ symmetric matrix $(g_{ij}(x))_{i, j = 0}^3$. In particular, we have 
\beq
g(x) = \sum_{0\leq \alpha<\beta\leq 3} g_{\alpha\beta}(x)  F^{\alpha\beta} + \sum_{\alpha = 0}^3 g_{\alpha\alpha}(x) G^\alpha.
\eeq
There is no $E^{\alpha\beta}$ term because $g$ is symmetric. Since $g$ is non-degenerate, at least one of the coefficients of $F^{\alpha\beta}, G^\alpha$ is non-zero. Without loss of generality, we assume that $g_{00}(x_0)\neq 0$ for $x_0\in M$. The other cases can be dealt in the same way. Then there is a neighborhood $U$ of $x_0$ so that $g_{00}(x)$ is non-vanishing there. Thus, for any $x\in U$, we can use $\{g(x),  E^{\alpha\beta}, F^{\alpha \beta}, G^1, G^2, G^3\}$ as a basis for $\mcb$. Therefore, we can write quadratic form $w$ as 
\beqq\label{eqwx}
w(x) = \sum_{0\leq \alpha<\beta\leq 3}a_{\alpha\beta}(x)E^{\alpha\beta} + \sum_{0\leq \alpha<\beta\leq 3}b_{\alpha\beta}(x)F^{\alpha\beta} + \sum_{\alpha = 1}^3 c_\alpha(x) G^\alpha + c_0(x)g(x).
\eeqq
We will show that for null forms, the coefficients of $F^{\alpha\beta}, G^1, G^2, G^3$ are all zero. 

Notice that for any $\xi \in C^\infty(M; TM)$, we have
\beq
F^{\alpha \beta}(\xi, \xi) = 2 \xi^\alpha\xi^\beta, \ \ G^\alpha(\xi, \xi) = (\xi^\alpha)^2, \ \ 0\leq \alpha, \beta \leq 3.
\eeq
For $\xi$ as a row vector, we define a symmetric matrix $A = \xi \xi^T$ i.e.\ $A_{ij} = \xi^i\xi^j, 0\leq i\leq j\leq 3$. Therefore, from $g(\xi, \xi) = 0$, we obtain that 
\beqq\label{eqga1}
 \mcf_1\doteq \sum_{0\leq \alpha<\beta\leq 3} 2 g_{\alpha\beta}(x) A_{\alpha\beta}  + \sum_{\alpha = 0}^3 g_{\alpha\alpha}(x)A_{\alpha\alpha} = 0.
\eeqq
For any fixed $x\in U$, the solution set of $\mcf_1 = 0$, denoted by $\mcs$, is a $9$-dimensional subspace of symmetric matrices. If $w$ satisfies the null condition \eqref{cond}, from \eqref{eqwx}, we should have 
\beq
\mcf_2\doteq \sum_{0\leq \alpha<\beta\leq 3}b_{\alpha\beta}(x)A_{\alpha\beta} + \sum_{\alpha = 1}^3 c_\alpha(x) A_{\alpha\alpha} = 0.
\eeq
Notice that this equation is linearly independent to \eqref{eqga1} as $g_{00}(x)\neq 0$, unless $b_\bullet=c_\bullet = 0$. Thus we conclude that $\mcf_2 = 0$ on $\mcs$ only if $b_\bullet = c_\bullet =0$ at $x$. Since this is true for any $x \in U$, we obtain that 
\beq
w(x) = \sum_{0\leq \alpha<\beta\leq 3}a_{\alpha\beta}(x)E^{\alpha\beta}  + c_0(x)g(x). 
\eeq
 This finishes the proof. 
\epf

We remark that using this lemma, we can locally write the quadratic form $w$ satsfying condition \eqref{eqw} as
\beqq\label{eqwnew}
w(x, u, \nabla_g u) =  C_0(x)g(x, \nabla_g u) + u C_1(x)g(x, \nabla_g u) +  u^2 \mcm(x,  \nabla_g u) + o(|u|^4 \cdot |\nabla_g u|^4).
\eeqq
However, this would not be true if $u$ is vector valued.  

Next, we discuss the well-posedness of \eqref{eqsem}. For the Cauchy problem, the local and global  existence of solutions to \eqref{eqsem} are extensively studied, see e.g.\ \cite{So, Tao}. For the source problem \eqref{eqsem}, we shall apply the results in \cite[Section 3.1.2]{KLU}  and \cite[Appendix B]{KLU1}, which are proved for more general second order quasilinear systems. These regularity  results may not be optimal but are sufficient for our purpose.

Let $B\subset \mbr^3$ be a compact set, $T_0>0$ and $f\in C^r([0,T_0], H^{s}(B))\cap C^{r+1}([0, T_0], H^{s-1}(B)), r\geq 0, s\geq 1$. If $m = r+ s\geq 4$ is even and $f$ is small enough, there is a unique solution $u\in C_0^r([0,T_0], H_0^{s}(\mbr^3))\cap C_0^{r+1}([0, T_0], H_0^{s-1}(\mbr^3))$ to \eqref{eqsem} such that  
\beqq\label{eqestimate}
\|u\|_{C_0^r([0,T_0], H_0^{s}(\mbr^3))\cap C_0^{r+1}([0, T_0], H_0^{s-1}(\mbr^3))} \leq C \|f\|_{C^r([0,T_0], H^{s}(B))\cap C^{r+1}([0, T_0], H^{s-1}(B))},
\eeqq
for some constant $C>0$ depending on $T_0$, see equation (27) of \cite{KLU1}. Hereafter, $C$ denotes a generic constant.  It is convenient to state the result in terms of Sobolev regularities. In particular,  from \eqref{eqestimate} we know that $u\in H^{m}(M(T_0))$ and $\|u\|_{H^{m}(M(T_0))}\leq C\|f\|_{H^{m}(M(T_0))}$ for $m\geq 4$ even.\\

Finally if $f$ depends on a small parameter $\eps,$ we find the asymptotic expansion of $u$ as $\eps \rightarrow 0.$ Let $f_i \in H_{\comp}^{8}(M), i = 1, 2, 3, 4$ be compactly supported in $M(T_0)\backslash M(0)$. Consider the linear wave equations
\beqq\label{eqlin}
\begin{gathered}
\square_g v_i = f_i \ \ \text{in } M(T_0), \\
v_i = 0 \text{ in } M(T_0)\backslash J_{g}^+(\supp(f_i)).
\end{gathered}
\eeqq
Since $g$ is globally hyperbolic, we know from for example \cite{Bar, Fr} that $\square_g$ has a causal inverse which we denote by $Q_g$, and $Q_g : H^m_{\comp}(M(T_0)) \rightarrow  H^{m+1}_{\loc}(M(T_0))$, see for example \cite[Prop.\ 5.6]{DUV} or \cite[Theorem 3.3]{GrU90}. So we have  $v_i = Q_g(f_i) \in H_{\loc}^{9}(M(T_0))$. Also, by the finite speed of propagation, we know that $v_i \in H^9_{\comp}(M(T_0))$ as $f_i$ are compactly supported. Let $\eps_i > 0, i = 1, 2, 3, 4$ be four small parameters and take
\beq
f = \sum_{i = 1}^4 \eps_i f_i \in H_{\comp}^{8}(M(T_0))\subset C_0^4(M(T_0))
\eeq
to be the source term in \eqref{eqsem}.  By the stability estimate, the solution $u$ to the nonlinear equation \eqref{eqsem} satisfies
\beq
\|u\|_{H^{8}(M(T_0))} \leq C(\eps_1+ \eps_2+\eps_3+\eps_4) \sum_{i = 1}^4\|f_i\|_{H^{8}(M(T_0))}.
\eeq
Also, the solution $v$ to the linearized equation \eqref{eqlin} with source $f$ is $v = \sum_{i = 1}^4 \eps_i v_i \in H^{9}(M(T_0))$. 

Now we derive the asymptotic expansion of $u$ as $\eps_i\rightarrow0$. Using \eqref{eqsem} and \eqref{eqlin}, we obtain that
\beq
\square_g (u - v) + w(x, u, \nabla_g u) = 0,
\eeq
which gives 
\beqq\label{eqite}
\begin{split}
u &= v - Q_g[w(x, u, \nabla_g u)]\\
& = v - Q_g[C_0(x)g(x, \nabla_g u) + C_1(x)u g(x, \nabla_g u) + u^2 \mcm(x, \nabla_g u)] + o(\eps_1\eps_2\eps_3\eps_4).
\end{split}
\eeqq
Here we used the expression \eqref{eqwnew}. Since $u\in H^8(M(T_0))$, we know that $\nabla_g u$ belongs to $H^7(M(T_0))$ which is an algebra. Thus the multiplication makes sense and $w(x, u, \nabla_g u)\in H^7(M(T_0))$. 
Therefore, the remainder terms $o(\eps_1\eps_2\eps_3\eps_4)$ in \eqref{eqite} are in $H^8(M(T_0)).$ Below, when it becomes necessary, we write the quadratic form $w$ as $w(x, u, \nabla_g u, \nabla_g u)$.
 
We look for the asymptotic expansion
\beq
u = v + \mcu^{(2)} + \mcu^{(3)} + \eps_1\eps_2\eps_3\eps_4\mcu^{(4)} + o(\eps_1\eps_2\eps_3\eps_4),
\eeq
where $\mcu^{(2)}$ denotes the collection of terms of the order $\eps_i \eps_j, i, j = 1, 2, 3, 4$ and $\mcu^{(3)}$ denotes the collection of terms of the order $ \eps_i\eps_j\eps_k, i, j, k = 1, 2, 3, 4$. We are interested in $\mcu^{(4)}$. Observe that this term can be determined from the source-to-solution map as  
\beqq\label{eqmcu4}
\mcu^{(4)} = \p_{\eps_1}\p_{\eps_2}\p_{\eps_3}\p_{\eps_4}u|_{\{\eps_1 = \eps_2 = \eps_3 = \eps_4 = 0\}} = \p_{\eps_1}\p_{\eps_2}\p_{\eps_3}\p_{\eps_4}L(f)|_{\{\eps_1 = \eps_2 = \eps_3 = \eps_4 = 0\}}.
\eeqq
To find the expression of $\mcu^{(4)}$, we substitute $u$ back to the right hand side of \eqref{eqite} and use $v = \sum_{i = 1}^4 \eps_i v_i$ to find order $\eps_1\eps_2\eps_3\eps_4$ terms. First, from the forth order term in \eqref{eqite}, we obtain
\beq
\mcm_1 =   -\sum_{(i, j, k, l) \in \sigma(4)}Q_g[ v_i v_j \mcm(x, \nabla_g v_k, \nabla_g v_l)],
\eeq
where $\sigma(4)$ denotes the set of all permutations of $(1, 2, 3, 4)$. 

Next, we consider  terms involving $C_1(x)$ i.e.\ the cubic term in \eqref{eqite}:
\beq
\begin{split}
\mcm_2 &= \sum_{(i, j, k, l) \in \sigma(4)} 2Q_g[C_0(x)g\big( \nabla_g v_i, \nabla_g Q_g\big( C_1(x)v_j g(\nabla_g v_k, \nabla_g v_l)\big)\big)] \\
& + \sum_{(i, j, k, l) \in \sigma(4)} Q_g[C_1(x)Q_g\big( C_0(x)g(\nabla_g v_i, \nabla_g v_j)\big) g(\nabla_g v_k, \nabla_g v_l)] \\
&+ \sum_{(i, j, k, l) \in \sigma(4)} 2Q_g[C_1(x)v_i g\big(\nabla_g v_j, \nabla_g Q_g(C_0(x)g(\nabla_g v_k, \nabla_g v_l))\big)]. 
\end{split}
\eeq
Finally, consider terms involving only $C_0(x)$ i.e.\ the quadratic term in \eqref{eqite}. To get these terms, we need to repeat the iteration twice using \eqref{eqite}.
\beq
\begin{gathered}
\mcm_3 = -\sum_{(i, j, k, l) \in \sigma(4)} 4Q_g[C_0(x)g(\nabla_g v_i, \nabla_g Q_g [C_0(x)g(\nabla_g v_j, \nabla_g Q_g g(\nabla_g v_k, \nabla_g v_l))])]\\
 -\sum_{(i, j, k, l) \in \sigma(4)}Q_g[C_0(x)g(\nabla_g Q_g \big( C_0(x) g(\nabla_g v_i, \nabla_g v_j)\big), \nabla_g Q_g \big(C_0(x)g(\nabla_g v_k, \nabla_g v_l)\big)].
\end{gathered}
\eeq
These are all the order $\eps_1\eps_2\eps_3\eps_4$ terms, hence we obtain that $\mcu^{(4)} = \mcm_1 + \mcm_2 + \mcm_3$. For convenience, we denote the $(i, j, k, l)$ term in $\mcu^{(4)}$ by $\mcu^{(4)}_{ijkl}$ and the term in $\mcm_\bullet, \bullet = 1, 2, 3$ by $\mcm_\bullet^{ijkl}$. We remark that these terms involve the multiplication of four linear waves $v_i$. We will analyze the singularities in $\mcu^{(4)}$ when $v_i$ has conormal type of singularities.

\section{The construction of distorted plane waves}\label{secdist}
We follow \cite{KLU} to construct distorted plane waves and before that, we recall some preliminaries. For $M = \mbr^{4},$ the cotangent bundle $T^*M$ is a symplectic manifold equipped with the canonical two form $w$. In local coordinates $(x, \xi)$ for $T^*M$, it is given by $\omega = d\xi\wedge dx$. A submanifold $\La \subset T^*M$ is called Lagrangian if the dimension of $\La$ is $4$ and $w$ vanishes on $\La$. Let $K\subset M$ be a submanifold. The conormal bundle of $K$ is defined as 
\beq
N^*K = \{(x, \xi)\in T^*M\backslash 0: \langle \xi, \theta \rangle = 0, \theta \in T_xK\},
\eeq
where $0$ represents the zero section of $T^*M$. It is known and easy to verify that $N^*K$ is a Lagrangian submanifold.  Now we briefly review the notion of (paired) Lagrangian distributions. Our references are \cite{DUV, Ho3, GU, GrU90, GrU93}.

Assume that $X$ is a $n$-dimensional smooth manifold and $\La$ is a smooth conic Lagrangian submanifold of $T^*X\backslash 0$.  We denote by $I^\mu(X; \La)$ the space of Lagrangian distributions of order $\mu$ on $X$ associated with $\La$, see e.g.\ \cite{Ho3, Ho4}. We also  abbreviate the notation by $I^\mu(X; \La) = I^\mu(\La)$ when the base manifold is not important. For example, if $Y$ is a submanifold of $X$, the conormal distributions to $Y$ are defined as $I^\mu(N^*Y)$.  Locally, Lagrangian distributions can be represented as oscillatory integrals. Let $\phi(x, \xi): U\times \mbr^d \rightarrow \mbr$ be a smooth non-degenerate phase function (homogeneous of degree $1$ in $\xi$) that  parametrizes $\La$ over an open set $U$ i.e. $\{(x, d_x\phi)\in T^*_UX\backslash 0: x\in U, d_\xi \phi = 0\} \subset \La.$ We can write $u\in I^{\mu}(\La)$ as a finite sum of oscillatory integrals
\beq
\int_{\mbr^N} e^{i\phi(x, \xi)} a(x, \xi) d\xi, \ \ a\in S^{\mu + \fnf - \frac{{d}}{2}}(U\times \mbr^{d}),
\eeq
where $S^\bullet(\bullet)$ denotes the standard symbol class, see \cite[Section 18.1]{Ho3}. For $u\in I^\mu(\La)$, we recall that the wave front set $\WF(u)\subset \La$ and $u\in H^s(X)$ for any $s< -\mu-\frac{n}{4}$, see \cite[Def.\ 25.1.1]{Ho4}. The distribution $u$ (as half densities on $M$) has a principal symbol $\sigma(u)$ defined invariantly on $\La$ (as half densities tensored with the Maslov factor), see \cite[Section 25.1]{Ho4}. {\em We emphasize that in our notation of principal symbols, we do not indicate the order but refer to the space of distributions for the order.} Also, for Lorentzian manifolds we consider, the density factor can be trivialized using the volume element $dg$. For conormal distributions, the Maslov factor in the symbol can be trivialized as well. 

Next consider paired Lagrangian distributions. Recall that two Lagrangians $\La_0, \La_1 \subset T^*X\backslash 0$ intersect cleanly at a codimension $k$ submanifold  if 
\beq
T_p\La_0\cap T_p\La_1  = T_p(\La_0\cap \La_1),\ \ \forall p\in \La_0\cap \La_1.
\eeq
The space of paired Lagrangian distribution associated with $(\La_0, \La_1)$ is denoted by $I^{p, l}(\La_0, \La_1)$. We shall not recall the definition here but only mention that if $u\in I^{p, l}(\La_0, \La_1)$, then $\WF(u)\subset \La_0\cup \La_1$. Also, away from the intersection $\La_0\cap \La_1$, $u \in I^{p+l}(\La_0\backslash \La_1)$ and $u \in I^p(\La_1\backslash \La_0)$ as Lagrangian distributions so $u$ has well-defined principal symbols $\sigma_{\La_0}(u)$ and $\sigma_{\La_1}(u)$ on corresponding Lagrangians (away from the intersections).  Our goal is to describe the Schwartz kernel of the causal inverse $Q_g$ as a paired Lagrangian distribution. {\em We remark that later we do not distinguish the notations of operators and their Schwartz kernels unless it is necessary.} 

Let $g^*$ be the dual Lorentzian metric on $T^*M$ induced by $g$. We let $\mcp(x, \xi) = |\xi|^2_{g^*}$ which is the  symbol of $\square_g$.  The characteristic set of $\mcp$ is defined as $\Sigma_g = \{(x, \xi)\in T^*M: \mcp(x, \xi) = 0\}$ and we notice that $\Sigma_g$ consists of light-like co-vectors. The Hamilton vector field $H_\mcp$  of $\mcp$ can be written in local coordinates as
\beq
H_\mcp = \sum_{i = 0}^3( \frac{\p \mcp}{\p \xi_i}\frac{\p }{\p x^i} - \frac{\p \mcp}{\p x^i}\frac{\p }{\p \xi_i}).
\eeq
The integral curves of $H_\mcp$ in $\Sigma_g$ are called null bicharacteristics and it is well-known that their projections to $M$ are geodesics. Now we move to the product manifold $M\times M$ and the cotangent bundle $T^*M\times T^*M$ as the Schwartz kernel of $Q_g$ is a distribution defined there. Let {$\pi_l: M\times M\rightarrow M$} be the 
projection to the left factor. We can regard $\mcp, \Sigma_g, H_\mcp$ as objects on product manifolds by pulling them back using {$\pi_l$}. Let $\diag = \{(z, z')\in M\times M: z = z'\}$ be the diagonal and 
\beq
N^*\diag = \{(z, \zeta, z', \zeta')\in T^*(M\times M)\backslash 0: z = z', \zeta' = -\zeta\}
\eeq 
be the conormal bundle of $\diag$. We define $\La_g$ to be the Lagrangian submanifold obtained by flowing out $N^*\diag\cap \Sigma_g$ under $H_\mcp$. It is proved in \cite{MU} (see also \cite{DUV}) that  $Q_g \in I^{-\frac{3}{2}, -\ha}(N^*\diag, \La_g)$.  

With these preparations, we construct distorted plane waves as in \cite{KLU}.  Let $L^+M$ be the bundle of future-pointing light-like vectors. For $(x^{(0)}, \theta^{(0)})\in L^+M$, we let $\gamma_{x^{(0)}, \theta^{(0)}}(t), t\geq 0$ be the unique geodesic from $x^{(0)}$ with direction $\theta^{(0)}$. Let $s_0, t_0>0$ be two small parameters (to be specified later). We define 
\beq
K(x^{(0)}, \theta^{(0)}; t_0, s_0) = \{\gamma_{x', \theta}(t)\in M(T_0); \theta \in \mco(s_0), t\in (0, \infty)\},
\eeq
where $(x', \theta') = (\gamma_{x^{(0)}, \theta^{(0)}}(t_0), \gamma_{x^{(0)}, \theta^{(0)}}'(t_0))$ and $\mco(s_0) \subset L^+_{x'}M$ is {an} open neighborhood of $\theta'$ consisting of ${\theta}\in L_{x'}^+M$ such that $\|{\theta} - \theta'\|_{g^+}< s_0$. 
It worths mentioning that as $s_0\rightarrow 0$, $K(x^{(0)}, \theta^{(0)}; t_0, s_0)$ tends to the geodesic $\gamma_{x^{(0)}, \theta^{(0)}}$. Now consider
\beq
Y(x^{(0)}, \theta^{(0)}; t_0, s_0) = K(x^{(0)}, \theta^{(0)}; t_0, s_0)  \cap \{t = 2t_0\},
\eeq
which intersects the geodesic at $\gamma_{x^{(0)}, \theta^{(0)}}(2t_0)$. We define $\La(x^{(0)}, \theta^{(0)}; t_0, s_0)$ to be the Lagrangian  obtained from flowing out $N^*K(x^{(0)}, \theta^{(0)}; t_0, s_0)\cap N^*Y(x^{(0)}, \theta^{(0)}; t_0, s_0)$ under the Hamilton vector field $H_\mcp$ in $\Sigma_g$. In the following, for simplicity we shall suppress the parameters $x^{(0)}, \theta^{(0)}, t_0, s_0$ in the notations  $K, Y$ and $\La$ etc.

For $\mu< -11$ , we take $f_0\in I^{\mu +1}(N^*Y)$   supported in a neighborhood $U$ of $\gamma_{x^{(0)}, \theta^{(0)}}\cap Y$ and that the principal symbol of $f_0$ vanishes outside of $\mco(s_0)\subset T^*M$. According to \cite[Lemma 3.1]{KLU}, $v_0 = Q_g f_0$ belongs to {$I^{\mu -\ha}(M\backslash Y; \La)$} and this is called the distorted plane wave. 
Moreover, the principal symbol of $v_0$ satisfies
\beqq\label{wavesym}
\sigma(v_0)(x, \xi) = \sigma_{\La_g}(Q_g)(x, \xi, y, \eta)\sigma(f_0)(y, \eta),
\eeqq
where $(x, \xi)$ and $(y, \eta)$ lie on the same bi-characteristics.  Since $\mu< -11$, we know that $f_0\in H^8_{\comp}(M)$ and $v_0\in H_{\loc}^9(M)$. So the regularity fits into the asymptotic analysis in Section \ref{asym}.

In this approach, we must pay attention to the conjugate points along $\gamma_{x^{(0)}, \theta^{(0)}}$, also called caustics. As discussed in \cite{KLU}, it is complicated to analyze the singularities past the caustics points. This difficulty will be overcome using another argument (where the parameter $t_0$ is useful). Let ${\tau_0}>0$ be such that $\gamma_{x^{(0)}, \theta^{(0)}}({\tau_0})$ is the first conjugate point of $x^{(0)}$ along $\gamma_{x^{(0)}, \theta^{(0)}}$. 
Then the exponential map $\exp_{x^{(0)}}$ is a local diffeomorphism from a neighborhood of $t\theta^{(0)}\in T_{x^{(0)}}M$ to a neighborhood of $\gamma_{x^{(0)}, \theta^{(0)}}(t)$ for $t< {\tau_0}$. Therefore, $K(x^{(0)}, \theta^{(0)}; t_0, s_0)$ is a codimension $1$ submanifold near $\gamma_{x^{(0)}, \theta^{(0)}}(t)$ and 
\beq
{\La (x^{(0)}, \theta^{(0)}; t_0, s_0) = N^*K(x^{(0)}, \theta^{(0)}; t_0, s_0) }
\quad \text{ near } \gamma_{x^{(0)}, \theta^{(0)}}(t) \text{ for }  t< {\tau_0}.
\eeq
In particular, before the first conjugate point of $x^{(0)}$ along $\gamma_{x^{(0)}, \theta^{(0)}}$, $v_0$ is a conormal distribution. For $s_0$ sufficiently small, the wave front of $v_0$ is concentrated near the geodesic.

\section{The nonlinear interaction of distorted plane waves}\label{secinter}
Assume that $x^{(j)} \in V$ and $(x^{(j)}, \theta^{(j)})\in L^+M, j = 1, 2, 3, 4$ are such that 
\beq
\gamma_{x^{(j)}, \theta^{(j)}}([0, t_0])\subset V, \ \ x^{(j)}(t_0)\notin J^+(x^{(k)}(t_0)) {\quad j\neq k},
\eeq
where $x^{(j)}(t) = \gamma_{x^{(j)}, \theta^{(j)}}(t)$. In particular, this means that the points $x^{(j)}(t_0)$ are causally independent.  We define $K_j = K(x^{(j)}, \theta^{(j)}; t_0, s_0)$ for $j = 1, 2, 3, 4$, and define $Y_j$ and $\La_j$ similarly as in Section \ref{secdist} (again, the parameters are suppressed for simplicity).   For $\mu <-11$, we let $f_j\in I^{\mu+1}(N^*Y_j)$ be constructed as $f_0$ in Section \ref{secdist} and $v_j = Q_g(f_j)\in I^{\mu-\ha}(M\backslash Y_j; \La_j)$. In the following, we analyze the singularities of $\mcu^{(4)}$ (defined in \eqref{eqmcu4}), where $v_i, i = 1, 2, 3, 4$ are distorted plane waves. In \cite[Section 3]{LUW}, the interaction terms of the form
\beqq\label{yterms}
\begin{gathered}
\mcy_1^{ijkl}  =Q_g(c(x)v_iv_jv_kv_l), \\
\mcy_2^{ijkl} = Q_g(b(x)v_iv_jQ_g(a(x)v_kv_l)),\ \ \mcy_3^{ijkl} = Q_g(a(x)v_iQ_g(b(x)v_jv_kv_l)),\\
\mcy_4^{ijkl} = Q_g(a(x)v_iQ_g(a(x)v_jQ_g(a(x)v_kv_l))), \ \ \mcy_5^{ijkl} = Q_g(a(x)Q_g(a(x)v_i v_j)Q_g(a(x)v_k v_l))
\end{gathered}
\eeqq
where $a(x), b(x), c(x)$ are smooth functions, are studied and the principal symbols are calculated. The only difference here is that we have derivatives $\nabla_g$ in the interaction term $\mcu^{(4)}$, for which we need to consider the action of $\nabla_g$ on conormal and paired Lagrangian distributions. We recall the lowering and raising of indices. The Lorentzian metric $g$ induces a natural isomorphism $T_xM\simeq T_x^*M$. For $\xi = (\xi_j)_{j = 0}^3 \in T_x^*M$, we let $\xi^\# = (\sum_{j = 0}^3g^{ij}\xi_j)_{i = 0}^3 \in T_xM$. For  $\xi = (\xi^j)_{j = 0}^3 \in T_xM$, we let $\xi^\flat = (\sum_{j = 0}^3g_{ij}\xi^j)_{i = 0}^3 \in T_x^*M$.

\begin{lemma}\label{dsym}
Let $\imath$ denote the imaginary unit i.e.\ $\imath^2 = -1.$
\begin{enumerate}
\item Let $Y$ be a codimension $1$ submanifold of $M$. For $u \in I^\mu(N^*Y)$, we have $\nabla_g u\in I^{\mu+1}(N^*Y)$. The principal symbol is given by ($i = 0, 1, 2, 3$)
\beq
\sigma((\nabla_g u)^i)(x, \xi) = \imath \xi^{\#, i}  \sigma(u)(x, \xi),  \ \ (x, \xi)\in N^*Y.
\eeq

\item 
Let ${Q}\in I^{p, l}(M\times M; N^*\diag, \La_g)$ and $\nabla_g$ act on the left factor of $M\times M$. We have ${\nabla_g} Q \in I^{p+1, l}(N^*\diag, \La_g).$
As a result, we have $\nabla_g Q \in I^{p+1+ l}(N^*\diag\backslash \La_g)$ and $\nabla_g Q \in I^{p+1}(\La_g\backslash N^*\diag)$. Moreover, their principal symbols are given by ($i = 0, 1, 2, 3$)
\beq
\begin{gathered}
\sigma((\nabla_g Q)^i)(x, \xi, x, -\xi) = \imath  \xi^{\#, i}  \sigma(Q)(x, \xi, x, -\xi), \text{ for } (x, \xi, x, -\xi) \in N^*\diag\backslash \La_g, \\
\sigma((\nabla_g Q)^i)(x, \xi, y, \eta) = \imath \xi^{\#, i} \sigma(Q)(x, \xi, y, \eta), \text{ for } (x, \xi, y, \eta) \in \La_g\backslash N^*\diag.
\end{gathered}
\eeq
\end{enumerate}
\end{lemma}
\bpf
(1) Locally near $Y$, we can choose local coordinates $x = (x', x''), x' \in\mbr^d, d = 1$ such that  $Y=\{x' = 0\}$. Then $N^*Y = \{x' = 0, \xi'' = 0\}$ where $\xi = (\xi', \xi'')$ denotes the dual variable. We can write $u(x)\in I^\mu(N^*Y)$ as an oscillatory integral
\beq
u(x) = \int_{\mbr^d} e^{\imath x'\cdot \xi'} a(x; \xi') d\xi', \ \ a\in S^{\mu + \fnf - \frac{d}{2}}(M\times (\mbr^d\backslash 0)), \ \ n = 4.
\eeq
Since $\nabla_g^i  = \sum_{i, j = 0}^{3} g^{ij}\p_{j}$, we obtain that
\beqq\label{eqngu}
\nabla^i_g u(x) = \int_{\mbr^d} e^{\imath x'\cdot \xi'}  [\imath \sum_{j = 0}^d g^{ij} \xi_j a(x; \xi')] d\xi' + \int_{\mbr^d} e^{\imath x'\cdot \xi'} \nabla^i_g a(x; \xi') d\xi'.
\eeqq
Since the amplitude in the first integral on the right hand side belongs to $S^{\mu + 1 + \fnf - \frac d 2}(M\times (\mbr^d\backslash 0))$ and $\nabla_g^i a \in S^{\mu + \fnf - \frac d 2}(M\times (\mbr^d\backslash 0)),$ we conclude that $\nabla_g u \in I^{\mu + 1}(N^*Y)$. The principal symbol can be read from  \eqref{eqngu}.

(2) The proof is similar to that of part (1) once we find the local oscillatory integral representation of $Q$. In fact, we will use the representation in \cite[Prop.\ 2.1]{GrU93}, see also \cite{MU}. We let $(x, y)$ be the local coordinates on $M\times M$ and $(\xi, \eta)$ be the dual variables. We can choose local coordinate $x = (x_1, x'), y= (y_1, y'), \xi = (\xi_1, \xi'), \eta = (\eta_1, \eta')$ such that the two intersecting Lagrangians are represented as
\beq
N^*\diag = \{x = y, \xi = -\eta\}, \ \ \La_g = \{x' = y', \xi' = -\eta', \xi_1 = \eta_1 = 0\}.
\eeq
Then we can write $Q$ as
\beq
Q(x, y) = \int_{\mbr^4} e^{\imath (x-y) \cdot \xi} b(x, y, \xi'; \xi_1) d\xi, \ \ b\in S^{p+\ha, l-\ha}(\mbr^{4+ 4}\times(\mbr^3\backslash 0) \times \mbr).
\eeq
Here $S^{*,*}(\bullet)$ denotes the product symbol space, see \cite{GrU93} for more details. On $N^*\diag\backslash \La_g$ i.e.\ $\xi_1\neq 0$, we actually have $b\in S^{p+l}(\mbr^{4+ 4}\times(\mbr^4\backslash 0))$ in the standard symbol space. On the other hand, on $\La_g\backslash N^*\diag$, we have that 
\beq
Q(x, y) = \int_{\mbr^3} e^{\imath (x'-y') \cdot \xi'} c(x, y, \xi') d\xi', 
\eeq
where for $x_1\neq y_1$, 
\beq
c(x, y, \xi') = \int_{\mbr} e^{\imath (x_1 -y_1)\xi_1} b(x, y; \xi', \xi_1)d\xi_1 \in S^{p+\ha} (\mbr^{4+4}\times (\mbr^3\backslash 0)).
\eeq
Now we can compute $\nabla_g Q$ as in part (1) and find its symbols. This completes the proof of the lemma. 
\epf

To describe the singularities in $\mcu^{(4)}$ produced by the nonlinear interaction of four distorted plane waves, we need to consider the following two issues. First, we know from \cite{KLU, LUW} that the interaction of three distorted plane waves could produce conic type singularities which are not used in solving the inverse problem. Let $
\La^{(3)} = \cup_{1\leq i<j<k\leq 4} N^*(K_i\cap K_j\cap K_k), 
$
and $\La^{(3), g}$ be the flow out of $\La^{(3)}\cap \Sigma_g$ under the Hamiltonian flow. Then the (new) singularities due to the interaction of three distorted plane waves are contained in $\La^{(3)}\cup \La^{(3), g}$. To include the singularities on $\La_i$, we denote 
$
\Theta = (\cup_{i = 1}^4\La_i) \cup \La^{(3)}\cup \La^{(3), g}.
$
We let $\pi: T^*M\rightarrow M$ be the standard projection and denote $\mck = \pi(\Theta)$, which is a subset of $M$ and contains the singular support of the singularities in $\mcu^{(4)}$ due to at most three wave interactions. Recall that the set $\Theta$ (hence $\mck$) by definition depends on the parameter $s_0$ and as $s_0\rightarrow 0$, it tends to a set of Hausdorff dimension $2$. Eventually this set become relatively small compared to the singular support of $\mcu^{(4)}$. 

The second issue is the conjugate points. Let ${\tau_j}, j = 1, 2, 3, 4$ be such that $\gamma_{x^{(j)}, \theta^{(j)}}({\tau_j})$ is the first conjugate point of $x^{(j)}$ along the geodesics and ${\tau_{\textrm{min}}} = \min_{j = 1, 2, 3, 4}({\tau_j})$. We see that before $\gamma_{x^{(j)}, \theta^{(j)}}({\tau_{\textrm{min}}})$, we have
$v_j \in I^{\mu-\ha}(M\backslash Y_j; N^*K_j).$ To avoid the complexities beyond the first conjugate points, we consider the interactions only in the following set  
\beq
\begin{gathered}
\mathcal{N}((\vec x, \vec \theta), t_0) = M(T_0)\backslash \cup_{j = 1}^4 J^+(\gamma_{x^{(j)}, \theta^{(j)}}({\tau_j})), 
\end{gathered}
\eeq
i.e.\ away from the causal future of the conjugate points, where $\vec x = (x^{(1)}, x^{(2)}, x^{(3)}, x^{(4)}),  \vec \theta = (\theta^{(1)}, \theta^{(2)}, \theta^{(3)}, \theta^{(4)}).$

Now we state the main result  about the singularities of $\mcu^{(4)}$. The proof (omitted here) is the same as that of \cite[Prop.\ 4.2]{LUW} by adjusting the orders.
\begin{prop}\label{inter5}
Let  $v_i \in I^{\mu-\ha}(N^*K_i), \mu < -11, i = 1, 2, 3, 4$ be the distorted plane waves constructed in the beginning of this section, and let $\mcu^{(4)}$ be the interaction term defined using $v_i$. For $q_0\in M$  we let $\La_{q_0} = T_{q_0}^*M\backslash 0$ and $\La^g_{q_0}$ be its flow out.  For $s_0>0$ sufficiently small, we have
\begin{enumerate}
\item If $\cap_{j = 1}^4\gamma_{x^{(j)}, \theta^{(j)}}(t) = \emptyset$ for $t< {\tau_{\textrm{min}}}$,  then $\mcu^{(4)}$ is smooth in $\mathcal{N}((\vec x, \vec \theta), t_0)\backslash \mck$;
\item If $\cap_{j = 1}^4\gamma_{x^{(j)}, \theta^{(j)}}(t) = q_0$ for $t< {\tau_{\textrm{min}}}$ and the corresponding tangent vectors at $q_0$ are linearly independent, then in $\mathcal{N}((\vec x, \vec \theta), t_0)\backslash \mck$, we have $\mcu^{(4)} \in I^{4\mu + \frac{7}{2}}(\La^g_{q_0}\backslash\Theta)$. 
\end{enumerate}
\end{prop}

Finally, we show that the leading singularities of $\mcu^{(4)}$ is not always vanishing. 
\begin{prop}\label{4sym}
Consider the setting of Prop.\ref{inter5} (2). Suppose the nonlinear term $w$ of \eqref{eqsem} is written in local coordinate in the form \eqref{eqwnew}. Let $(q, \eta)\in \La^g_{q_0}\backslash \Theta$ be joined to $(q_0, \zeta), \eta\in T_{q_0}^*M$ by bicharacteristics. We can write $\zeta = \sum_{i = 1}^4 \zeta^{(i)}$ for $\zeta^{(i)} \in N^*_{q_0}K_i$   linearly independent. Then the principal symbol of $\mcu^{(4)}$ in $I^{4\mu + \frac{7}{2}}(\La^g_{q_0}\backslash\Theta)$ at $(q, \eta)$ is given by 
\beq
\begin{gathered}
  \sigma_{\La_{q_0}^g}(\mcu^{(4)})(q, \eta) =   -(2\pi)^{-3} \sigma(Q_g)(q, \eta, q_0, \zeta) \mcp( q_0, \zeta^{(1)}, \zeta^{(2)},  \zeta^{(3)},  \zeta^{(4)}) \prod_{i = 1}^4\sigma(v_i)(q_0, \zeta^{(i)}),\\
  \text{where } \mcp  = 2 \big(\mcm(q_0, \zeta^\#, \zeta^\#) -\sum_{i = 1}^4 \mcm(q_0, \zeta^{(i), \#}, \zeta^{(i), \#}) \big).
  \end{gathered}
\eeq
Moreover, the function $\mcp(q_0, \cdot)$ is  non-vanishing on any open set of  
\beq
\begin{gathered}
\mcx = \{( \zeta^{(1)},  \zeta^{(2)},   \zeta^{(3)},  \zeta^{(4)})\in (L_{q_0}^{*}M\backslash\{0\})^4:  \zeta^{(i)} \text{ are linearly independent} 
 \text{ and $\sum_{i = 1}^4 \zeta^{(i)}$ is light-like}\},
\end{gathered}
\eeq
where $L_{q_0}^*M$ denotes the set of light-like co-vectors at $q_0$. 
\end{prop}
 
\bpf 
First of all, we compute the principal symbol of $\mcu^{(4)}$.  In Section 3.5 of \cite{LUW}, the principal symbols of the terms $\mcy_\bullet^{ijkl}, \bullet = 1, 2, 3, 4, 5$ in \eqref{yterms} are found explicitly, and we recall them here. Consider the symbols at $(q, \eta)\in \La_{q_0}^{g}\backslash \Theta$ which is joined with $(q_0, \zeta)\in \La_{q_0}$ by bi-characteristics. We can write $\zeta = \sum_{i = 1}^4\zeta^{(i)}$ where $\zeta^{(i)}\in N^*_{q_0}K_i$. Let $A_i$ be the principal symbols of $v_i$. By Prop.\ 3.12 of \cite{LUW}, we get 
\beq
 \sigma_{\La_{q_0}}(\mcy_\bullet^{ijkl})(q, \eta) =  (2\pi)^{-3} \sigma_{\La_{g}}(Q_{g})(q, \eta, q_0, \zeta)\mcp_\bullet(\zeta^{(1)}, \zeta^{(2)}, \zeta^{(3)}, \zeta^{(4)}) \prod_{i = 1}^4A_i(q_0, \zeta^{(i)}), \ \ \bullet = 1, 2, 3, 4, 5, 
\eeq
where $ \mcp_1 = c(q_0)$ and 
 \beq
 \left. \begin{array}{ll}
 \mcp_2 =   a(q_0)b(q_0)   |\zeta^{(k)}+\zeta^{(l)}|_{g^*(q_0)}^{-2}, & \mcp_3 = a(q_0)  b(q_0) |\zeta^{(j)}+\zeta^{(k)}+\zeta^{(l)}|_{g^*(q_0)}^{-2}, \\
  \mcp_4 =   a^3(q_0)   |\zeta^{(j)}+\zeta^{(k)}+\zeta^{(l)}|_{g^*(q_0)}^{-2}  |\zeta^{(k)}+\zeta^{(l)}|_{g^*(q_0)}^{-2}, & \mcp_5  =  a^3(q_0) |\zeta^{(k)}+\zeta^{(l)}|_{g^*(q_0)}^{-2}   |\zeta^{(i)}+\zeta^{(j)}|_{g^*(q_0)}^{-2}. 
  \end{array}\right.
 \eeq 

To find the principal symbols of $\mcu^{(4)} = \mcm_1 + \mcm_2 + \mcm_3$, we just need to use Lemma \ref{dsym}  to take into account the derivatives. We start from $\mcm_1$ in which the terms are like $\mcy_1^\bullet$. We have
\beq
\begin{gathered}
\sigma(\mcm_1)  = -(2\pi)^{-3} \sigma_{\La_g}(Q_g)(q, \eta, q_0, \zeta) \cdot \mcp \cdot \prod_{i = 1}^4A_i(q_0, \zeta^{(i)}),\\
\text{where } \mcp  = \sum_{(i, j, k, l) \in \sigma(4)} \mcm(q_0, \zeta^{(k), \#}, \zeta^{(l), \#}) =2 \big(\mcm(q_0, \zeta^\#, \zeta^\#) -\sum_{i = 1}^4 \mcm(q_0, \zeta^{(i), \#}, \zeta^{(i), \#}) \big).
\end{gathered} 
\eeq

Next, consider $\mcm_2$ in which the terms are like $\mcy_2^\bullet, \mcy_3^\bullet$. We get
\beq
 \begin{split}
&\sigma(\mcm_2)(q, \eta)   =  -(2\pi)^{-3}  C_0(q_0)C_1(q_0) \sigma_{\La_g}(Q_g)(q, \eta, q_0, \zeta)\cdot \mca \cdot \prod_{i = 1}^4A_i(q_0, \zeta^{(i)})\\
\text{where }  \mca  &=  \sum_{(i, j, k, l)\in \sigma(4)} [2\frac{g(\zeta^{(i)}, \zeta^{(j)}+\zeta^{(k)}+\zeta^{(l)})}{|\zeta^{(j)}+\zeta^{(k)}+\zeta^{(l)}|_{g^*(q_0)}^{2}} \cdot  g(\zeta^{(k)}, \zeta^{(l)}) + \frac{g(\zeta^{(i)}, \zeta^{(j)}) }{ |\zeta^{(i)}+\zeta^{(j)}|_{g^*(q_0)}^{2}}g(\zeta^{(k)}, \zeta^{(l)}) \\
 & \phantom{aaaaaaaaa} +2\frac{g(\zeta^{(k)}, \zeta^{(l)})}{|\zeta^{(k)}+\zeta^{(l)}|_{g^*(q_0)}^{2}} \cdot  g(\zeta^{(j)}, \zeta^{(k)}+ \zeta^{(l)}) ]\\
 &=   \sum_{(i, j, k, l)\in \sigma(4)} [2\frac{g(\zeta^{(i)}, \zeta^{(j)}+\zeta^{(k)}+\zeta^{(l)})}{|\zeta^{(j)}+\zeta^{(k)}+\zeta^{(l)}|_{g^*(q_0)}^{2}} \cdot  g(\zeta^{(k)}, \zeta^{(l)}) + \ha g(\zeta^{(k)}, \zeta^{(l)}) + g(\zeta^{(j)}, \zeta^{(k)}+ \zeta^{(l)}) ].
 \end{split}
\eeq
Here we abused the notations that the vectors $\zeta^{(i)}$ inside $g$ should be regarded as tangent vectors while the vectors in $|\cdot|_{g^*(q_0)}$ are the cotangent vectors. But we have for $\xi, \eta \in T_{q_0}M$ that
$
g(\xi, \eta) = \sum_{i, j = 0}^3 g_{ij} \xi^i \eta^j = \sum_{i, j = 0}^3 g^{ij}\xi_i \eta_j = g^*(\xi^\flat, \eta^\flat).
$

Consider the first term in the summation of $\mca$. For $i = 1$, we sum in $(j, k, l) $ over $\sigma(3)$ the permutations of $(2, 3, 4)$. We have 
\beq
\begin{split}
 \sum_ {(j, k, l) \in \sigma(3)}  \frac{g^*(\zeta^{(k)}, \zeta^{(l)})}{|\zeta^{(j)}+\zeta^{(k)}+\zeta^{(l)}|^2_{g^*(q_0)}} =   \frac{g^*(\zeta^{(2)}, \zeta^{(3)}) + g^*(\zeta^{(2)}, \zeta^{(4)}) +  g^*(\zeta^{(3)}, \zeta^{(4)})}{g^*(\zeta^{(2)}, \zeta^{(3)}) + g^*(\zeta^{(2)}, \zeta^{(4)}) +  g^*(\zeta^{(3)}, \zeta^{(4)})} = 1.
 \end{split}
 \eeq
Therefore, for $i = 1$, the summation of the first term in $\mca$ is  $2g(\zeta^{(1)}, \zeta^{(2)}+\zeta^{(3)}+\zeta^{(4)})$. Similarly, we can compute for $i = 2, 3, 4$. So the summation of the first term in $\mca$ gives
\beq
\begin{split}
 &\sum_{(i, j, k, l)\in \sigma(4)}  2\frac{g(\zeta^{(i)}, \zeta^{(j)}+\zeta^{(k)}+\zeta^{(l)})}{|\zeta^{(j)}+\zeta^{(k)}+\zeta^{(l)}|_{g^*(q_0)}^{2}} \cdot  g(\zeta^{(k)}, \zeta^{(l)}) \\
& =  4[g^*(\zeta^{(1)}, \zeta^{(2)}+\zeta^{(3)}+\zeta^{(4)}) + g^*(\zeta^{(2)}, \zeta^{(1)}+\zeta^{(3)}+\zeta^{(4)}) + g^*(\zeta^{(3)}, \zeta^{(1)}+\zeta^{(2)}+\zeta^{(4)})\\
  &\phantom{aaa}+ g^*(\zeta^{(4)}, \zeta^{(1)}+\zeta^{(2)}+\zeta^{(3)})] \\
& =  8[g^*(\zeta^{(1)}, \zeta^{(2)}) +g^*(\zeta^{(1)}, \zeta^{(3)}) + g^*(\zeta^{(1)}, \zeta^{(4)}) + g^*(\zeta^{(2)}, \zeta^{(3)}) + g^*(\zeta^{(2)}, \zeta^{(4)}) + g^*(\zeta^{(3)}, \zeta^{(4)})]\\
 & = 4 g^*(\zeta, \zeta) = 0
\end{split}
\eeq
Thus, we get
\[
\mca = \sum_{(i, j, k, l) \in \sigma(4)}\frac{5}{2} g^*(\zeta^{(k)}, \zeta^{(l)}) = 5 g^*(\zeta, \zeta) = 0,
\]
 hence the principal symbol $\sigma(\mcm_2)$ vanishes. 

Finally, consider the term $\mcm_3$ in which the terms are similar to $\mcy_4^\bullet, \mcy_5^\bullet$. We find that
\beq
\begin{split}
&\sigma(\mcm_3)(q, \eta)    = (2\pi)^{-3} [C_0(q_0)]^3\sigma_{\La_g}(Q_g)(q, \eta, q_0, \zeta)\cdot \mcb \cdot \prod_{i = 1}^4A_i(q_0, \zeta^{(i)}),\\
\text{where } \mcb = & \sum_{(i, j, k, l)\in \sigma(4)} [4\frac{g(\zeta^{(i)}, \zeta^{(j)}+\zeta^{(k)}+\zeta^{(l)})}{|\zeta^{(j)}+\zeta^{(k)}+\zeta^{(l)}|_{g^*(q_0)}^{2}} \cdot \frac{g(\zeta^{(j)}, \zeta^{(k)}+ \zeta^{(l)})}{|\zeta^{(k)}+\zeta^{(l)}|_{g^*(q_0)}^{2}} g(\zeta^{(k)}, \zeta^{(l)}) \\
 & \phantom{aaaaaa}+ \frac{g(\zeta^{(k)}+\zeta^{(l)}, \zeta^{(i)}+\zeta^{(j)})}{|\zeta^{(k)}+\zeta^{(l)}|_{g^*(q_0)}^{2} |\zeta^{(i)}+\zeta^{(j)}|_{g^*(q_0)}^{2}}g(\zeta^{(i)}, \zeta^{(j)})g(\zeta^{(k)}, \zeta^{(l)})]\\
 = & \sum_{(i, j, k, l)\in \sigma(4)}  [2\frac{g^*(\zeta^{(i)}, \zeta^{(j)}+\zeta^{(k)}+\zeta^{(l)})}{|\zeta^{(j)}+\zeta^{(k)}+\zeta^{(l)}|^2_{g^*(q_0)}}g^*(\zeta^{(j)}, \zeta^{(k)}+\zeta^{(l)})  + \frac{1}{4}g^*(\zeta^{(k)}+\zeta^{(l)}, \zeta^{(i)}+\zeta^{(j)})].
\end{split}
\eeq
We consider the summation of the first term.  Observe that for fixed $i = 1$ and $(j, k, l) \in \sigma(3)$ the set of permutations of $(2, 3, 4)$, we have
\beq
\begin{split}
 \sum_ {(j, k, l) \in \sigma(3)}  \frac{g^*(\zeta^{(j)}, \zeta^{(k)}+\zeta^{(l)})}{|\zeta^{(j)}+\zeta^{(k)}+\zeta^{(l)}|^2_{g^*(q_0)}} 
  =  & \ha \cdot 4 \cdot \frac{g^*(\zeta^{(2)}, \zeta^{(3)}) + g^*(\zeta^{(2)}, \zeta^{(4)}) +  g^*(\zeta^{(3)}, \zeta^{(4)})}{g^*(\zeta^{(2)}, \zeta^{(3)}) + g^*(\zeta^{(2)}, \zeta^{(4)}) +  g^*(\zeta^{(3)}, \zeta^{(4)})} = 2.
 \end{split}
 \eeq
The situation is the same for $i = 2, 3, 4$. Therefore, we obtain that
\beq
\begin{split}
\mcb &=   4[g^*(\zeta^{(1)}, \zeta^{(2)}+\zeta^{(3)}+\zeta^{(4)}) + g^*(\zeta^{(2)}, \zeta^{(1)}+\zeta^{(3)}+\zeta^{(4)}) + g^*(\zeta^{(3)}, \zeta^{(1)}+\zeta^{(2)}+\zeta^{(4)})\\
  &\phantom{aa}+ g^*(\zeta^{(4)}, \zeta^{(1)}+\zeta^{(2)}+\zeta^{(3)})]   + 2[g^*(\zeta^{(1)}+\zeta^{(2)}, \zeta^{(3)}+\zeta^{(4)}) + g^*(\zeta^{(1)}+\zeta^{(3)}, \zeta^{(2)}+\zeta^{(4)}) \\
  & \phantom{aa}+ g^*(\zeta^{(1)}+\zeta^{(4)}, \zeta^{(2)}+\zeta^{(3)})]\\
 &=  12[g^*(\zeta^{(1)}, \zeta^{(2)}) +g^*(\zeta^{(1)}, \zeta^{(3)}) + g^*(\zeta^{(1)}, \zeta^{(4)}) + g^*(\zeta^{(2)}, \zeta^{(3)}) + g^*(\zeta^{(2)}, \zeta^{(4)}) + g^*(\zeta^{(3)}, \zeta^{(4)})]\\
  &= 3 g^*(\zeta, \zeta) = 0
\end{split}
\eeq
Thus the principal symbol of $\mcm_3$ also vanishes.  Since $\mcu^{(4)} = \mcm_1 + \mcm_2 + \mcm_3$, we have $\sigma(\mcu^{(4)})(q, \eta) = \sigma(\mcm_1)(q, \eta)$ and we proved the first claim of the proposition.

It remains to show that $\mcp$  does not vanish on any open subset of $\mcx$. We start from
\beq
\tilde \mcx = \{( \zeta^{(1)},  \zeta^{(2)},   \zeta^{(3)},  \zeta^{(4)})\in (L_{q_0}^{*}M\backslash\{0\})^4:  \zeta^{(i)} \text{ are linearly independent}\}.
\eeq
This is an open subset of $(L_{q_0}^*M)^4$ which is a manifold of dimension $12$. We consider a smooth map $\mcf: \tilde \mcx \rightarrow \mbr$ given by $\mcf =  g(\zeta, \zeta),   \zeta = \sum_{i = 1}^4 \zeta^{(i)}.$ 
Then $\mcx = \mcf^{-1}(0)$. We compute the Jacobian of $\mcf$. Observe that 
\beq
\frac{\p}{\p \zeta^{(1)}}g(\zeta,\zeta) = (2g^{0j}\zeta_j, 2g^{1j}\zeta_j, 2g^{2j}\zeta_j, 2g^{3j}\zeta_j) = 2G\zeta,
\eeq
where  $G = (g^{ij}(q_0))$ is a $4\times 4$ matrix. Therefore, the Jacobian is
\beq
D\mcf = (    2G\zeta , 2G\zeta, 2G\zeta, 2G\zeta).
\eeq
Since $\zeta^{(i)}$ are linearly independent and $G$ is non-degenerate, we see that $D\mcf$ has constant rank $1$ so that $\mcx$ is a smooth manifold of dimension $11$, see e.g. \cite[Theorem 5.8]{Boo}. Now we consider the subset of $\mcx$ where $\mcp$ vanishes. Consider the map $(\mcf, \mcp) : \tilde \mcx  \rightarrow \mbr^2$. We compute the Jacobian of $\mcp$ as
\beq
D\mcp =4(GMG (\zeta - \zeta^{(1)}),  GMG(\zeta - \zeta^{(2)}), GMG(\zeta - \zeta^{(3)}), GMG (\zeta - \zeta^{(4)})), 
\eeq
where $M  = (\mcm^{ij}(q_0))$ is a $4\times 4$ matrix. If $D(\mcf, \mcp)$ has rank $1$, we can find constant $\beta$ such that 
\beq
\begin{gathered}
4 GMG (\zeta - \zeta^{(a)}) = 2 \beta   G\zeta, \ \ a = 1, 2, 3, 4.
\end{gathered}
\eeq
Summing over $a$, we get $3GMG \zeta =  2\beta G\zeta.$ Since we can choose $\zeta$ to be any light-like vector, we must have $M = \frac {2\beta}{3} G^{-1}$ at $q_0$. This in particular means that for $\zeta\in L_{q_0}M$, 
\beq
\mcm(\zeta, \zeta) = \zeta^t M \zeta =  \frac {2\beta}{3} \zeta^t G^{-1}\zeta  =  \frac {2\beta}{3} g(\zeta, \zeta) = 0.
\eeq
But we know by assumption \eqref{eqw} that $\mcm$ is not null. Therefore, the rank of $D(\mcf, \mcp)$ is $2$. As a result, $\mcp$ vanishes on a $10$ dimensional submanifold of $\mcx$, which means $\mcp$ cannot vanish on any open subset of $\mcx.$ This finishes the proof.
\epf

We remark that the proposition can be formulated as follows. Let $\zeta = \sum_{i = 1}^4 \zeta^{(i)} \in L^*_{q_0}M, \vec \zeta = (\zeta^{(i)})_{i = 1}^4\in \mcx$. For any neighborhood $W\subset L_{q_0}^*M$ of $\zeta$, it follows from Prop.\ \ref{4sym} that there is a neighborhood $\mcw \subset \mcx$ of $\vec \zeta$ such that $\mcp$ is not always vanishing on $\mcw$. Therefore, one can think $\mcp(q_0, \cdot)$ as a function defined on the set 
\beq
\begin{gathered}
\mcz =  \{( \zeta^{(1)},  \zeta^{(2)},   \zeta^{(3)},  \zeta^{(4)}, \zeta )\in (L_{q_0}^{*}M\backslash\{0\})^5:  \zeta^{(i)} \text{ are linearly independent} 
 \text{ and $\zeta = \sum_{i = 1}^4 \zeta^{(i)}$}\}.
\end{gathered}
\eeq
Then $\mcp(q_0, \cdot)$ is non-vanishing on any open subset of $\mcz$. This puts Prop.\ \ref{4sym} in the same form as \cite[Prop.\ 3.4]{KLU} or \cite[Prop.\ 3.4]{LUW}.

\section{Proof of the main theorems}\label{secinv}
We prove our main results essentially following the arguments in \cite{KLU} and \cite[Section 4]{LUW}. Since the proofs are very similar, we shall only go over the key components and refer the readers to the above works for details. 

To deal with conjugate points, we recall the earliest light observation set  introduced in \cite{KLU}. The light observation set of $q\in M$ in $V$ is defined as $\mcp_V(q) =(J^+(q)\backslash I^+(q))\cap V$ i.e.\ points in $V$ which are on the future pointing light-like curves from $q$. The earliest light observation set is defined as 
\beq
\begin{gathered}
\mce_V(q) = \{x\in \mcp_V(q): \text{there is no $y\in \mcp_V(q)$ and future-pointing time-like path} \\
 \text{$\alpha:[0, 1] \rightarrow V$ such that $\alpha(0) = y$ and $\alpha(1) = x$} \}\subset V,
 \end{gathered}
\eeq
see \cite[Def.\ 1.1]{KLU}. For $W\subset M$ open, we let $\mce_V(W) = \{\mce_V(q): q\in W\}.$ For exmaple, for the interaction point  $q_0$  in Prop.\ \ref{inter5}, we have $\mce_V(q_0) \subset \mathcal{N}((\vec x, \vec \theta), t_0). $ Also, the set $\mce_V(q_0)$ is not empty because $\mce_V(q_0)$ contains a $3$-dimensional submanifold as discussed in Section 2.2.1 of \cite{KLU}.  

\bpf[Proof of Theorem \ref{main}]
We first prove the determination of the conformal class. Our Prop.\ \ref{inter5} and Prop.\ \ref{4sym} are equivalent to Theorem 3.3 and Prop.\ 3.4 of \cite{KLU}.  By the arguments in Section 3.5 and Section 4 of \cite{KLU}, we can show that the source-to-solution map $L$ determines the earliest light observation set $\mce_V$ of a dense subset of $I(p_-, p_+)$ (here as in \cite{KLU}, we shall take the parameter $s_0\rightarrow 0$ in definition of $K_i$ so that $\pi(\Theta)$ tends to a set of Hausdorff dimension $2$). The problem is reduced to the inverse problem with passive measurements. It follows from Theorem 1.2 and Remark 2.2 of \cite{KLU}  that the differential structure of $I(p_-, p_+)$ and the conformal class of the metric can be uniquely determined up to diffeomorphisms. In the case when $g^{(i)}$ are Ricci flat, it follows from Corollary 1.3 of \cite{KLU2} that the conformal diffeomorphism is indeed an isometry.  

Finally, consider the case when $\mcm^{(1)}  = \mcm^{(2)}  = \mcm$ is a quadratic form independent of $x$. We already proved that the two metrics $g^{(i)}$ are conformal to each other. Without loss of generality, we can assume that $g^{(1)} = e^{2\gamma} g^{(2)}$. By linearizing the source-to-solution map, we can deduce from $L^{(1)}(f) = L^{(2)}(f)$ in $V$ that $g^{(1)}=  g^{(2)}$ in $V$ i.e. $\gamma = 0$ in $V$, see \cite[Remark 3.1]{KLU}. Now, assume $f = \sum_{i = 1}^4\eps_i f_i$ as constructed in Section \ref{secinter} and denote the interaction terms by
\beq
\mcu^{(4), \alpha} = \p_{\eps_1}\p_{\eps_2}\p_{\eps_3}\p_{\eps_4}L^{(\alpha)}(f)|_{\{\eps_1=\eps_2=\eps_3=\eps_4 =0\}}, \ \ \alpha = 1, 2.
\eeq
Then we know that $\mcu^{(4), 1}(q) = \mcu^{(4), 2}(q), q\in V.$ For any $q_0 \in I(p_-, p_+)$, we will compare the principal symbols of $\mcu^{(4), 1}$ and $\mcu^{(4), 2}$ on $\mce_{V}(q_0)\backslash \Theta$ for $s_0\rightarrow 0$ as in \cite[Section 4]{LUW}. We remark that since conformal transformations of Lorentzian metrics preserves light-like (pre)geodesics, the sets $\mce_V(q_0)$ are the same for $g^{(1)}, g^{(2)}$. We need the following result.

\begin{prop}[Prop.\ 4.5 of \cite{LUW}]\label{symconf}
Let $g, \tilde g$ be two Lorentzian metrics on $M$ such that $g = e^{2\gamma} \tilde g$ where $\gamma \in C^\infty(M).$ Let $Q_g, Q_{\tilde g}$ be the causal inverse of $\square_g, \square_{\tilde g}$ respectively. Then the Lagrangians $\La_g = \La_{\tilde g}$ and the principal symbols of $Q_g, Q_{\tilde g} \in I^{-2}(N^*\diag\backslash \La_g)$ satisfy $\sigma(Q_g) = e^{2\gamma} \sigma(Q_{\tilde g}).$ For their principal symbols in $I^{-\frac{3}{2}}(\La_g \backslash N^*\diag)$, we have
\beq
\sigma(Q_g)(x, \xi, y, \eta) = e^{-\gamma(x)} \sigma(Q_{\tilde g})(x, \xi, y, \eta) e^{3\gamma(y)},
\eeq
for $(x, \xi), (y, \eta)$ on the same bicharacteristics on $\La_g$.
\end{prop}

Now from Prop.\ \ref{4sym}, we have for any $(q, \eta)\in \La_{q_0}^g\backslash \Theta$ and $\alpha = 1, 2$ that
\beq
\begin{gathered}
  \sigma_{\La_{q_0}^g}(\mcu^{(4), \alpha})(q, \eta) = -(2\pi)^{-3} \sigma_{\La_g}(Q_g)(q, \eta, q_0, \zeta) \cdot \mcp^{(\alpha)} \cdot  \prod_{i = 1}^4A^{(\alpha)}_i(q_0, \zeta^{(i)}), \\
  \mcp^{(\alpha)}  = 2 \big(\mcm(\zeta^\#, \zeta^\#) -\sum_{i = 1}^4 \mcm(\zeta^{(i), \#}, \zeta^{(i), \#}) \big),
  \end{gathered}
\eeq
where the $\#$ operation is with respect to $g^{(\alpha)}$ and $A^{(\alpha)}_i$ are the principal symbols of $v_i^{(\alpha)} = Q_{g^{(\alpha)}}(f_i)$. These symbols satisfy
\beq
A^{(\alpha)}_i(q_0, \xi^{(i)}) = \sigma(Q_{g^{(\alpha)}})(q_0, \xi^{(i)}, x^{(i)}, \zeta^{(i)})B_i(x^{(i)}, \zeta^{(i)}), \ \ i = 1, 2, 3, 4, \alpha = 1, 2,
\eeq
where $x^{(i)}\in V$, $(q_0, \xi^{(i)})$  and $(x^{(i)}, \zeta^{(i)})$ are joined by bicharacteristics and $B_i$ are the principal symbols of $f_i$. Next, let's consider the conformal transformation of the symbols. Since $g^{(1)} = e^{2\gamma} g^{(2)}$, by Prop.\ \ref{symconf}, we have for $i = 1, 2, 3, 4$ that
\beqq\label{symrel}
\begin{gathered}
\sigma(Q_{g^{(1)}})(q, \eta, q_0, \zeta^{(i)}) = \sigma(Q_{g^{(2)}})(q, \eta, q_0, \zeta^{(i)})e^{3\gamma(q_0)},\\
\sigma(Q_{g^{(1)}})(q_0, \xi^{(i)}, x^{(i)}, \zeta^{(i)}) = e^{-\gamma(q_0)}\sigma(Q_{g^{(2)}})(q_0, \xi^{(i)}, x^{(i)}, \zeta^{(i)}).
\end{gathered}
\eeqq
Next, recall that for $\zeta\in T_{q_0}^*M$, $(\zeta^\#)^i = \sum_{j = 0}^3g^{(\alpha), ij}\zeta_j$. So we have $\mcp^{(1)} = e^{-4\gamma(q_0)}\mcp^{(2)}$. Finally, we obtain the following relation
\beqq\label{eqrelation}
 \sigma(\mcu^{(4), 1})(q, \eta) =  e^{- 5\gamma(q_0)} \sigma(\mcu^{(4), 2})(q, \eta), \ \ (q, \eta) \in \La_{q_0}^g\backslash \Theta.
\eeqq
According to Prop.\ \ref{inter5} and \ref{4sym}, we can choose $f_i$ so that the principal symbols $ \sigma(\mcu^{(4), \alpha})(q, \eta)$ are non-vanishing. Because $\mcu^{(4), 1} = \mcu^{(4), 2}$ on $V$, we know their principal symbols must be the same for $q\in V$. So we conclude that $e^{- 5\gamma(q_0)} = 1$. As this is true for all $q_0\in I(p_-, p_+)$, the proof is complete. 
\epf

\section*{Acknowledgment}
The authors would like to thank Prof.\ Gunther Uhlmann for suggesting the problem and for many helpful discussions. TZ was supported by NSF grant DMS-1501049 and Alfred P. Sloan Research Fellowship FR-2015-65641.


\end{document}